\documentclass[conference,a4paper]{IEEEtran}

\usepackage{comment}
\usepackage[utf8]{inputenc} 

\usepackage[T1]{fontenc}
\usepackage{url}              
\usepackage{cite}             

\usepackage[cmex10]{amsmath}  
\usepackage{mleftright}       
\mleftright                   

\usepackage{graphicx}         
\usepackage{booktabs}         
          
\usepackage{amssymb}

\newcommand{\cX}{{\mathcal{X}}}

\newcommand{\cY}{{\mathcal{Y}}}

\newcommand{\bone}{{\mathbf{1}}}

\usepackage[mathcal]{euscript}

\DeclareMathOperator*{\argmin}{arg\,min}

\DeclareMathOperator{\ev}{\mathbb{E}}
\DeclareMathOperator{\pr}{\mathbb{P}}

\newcounter{actr}
{\begin{list}{(\alph{actr})}{\usecounter{actr}}}{\end{list}}

\newcounter{ictr}
{\begin{list}{(\roman{ictr})}{\usecounter{ictr}}}{\end{list}}

{

}%
{

}

\newcommand{\abs}[1]{\left|#1\right|}

\newcommand{\defeq}{\triangleq}

\newcommand{\E}[1]{\ev\left[{#1}\right]}
\newcommand{\Ed}[2]{\ev_{#1}\left[{#2}\right]} 

\newcommand{\Prob}[1]{\pr\left({#1}\right)}
\newcommand{\prob}[1]{\Prob{#1}}

\newcommand{\reals}{\mathbb{R}}

\newcommand{\floor}[1]{\lfloor{#1}\rfloor}

\newcommand{\thmref}[1]{Theorem~\ref{#1}}

\newcommand{\lemref}[1]{Lemma~\ref{#1}}

\usepackage{todonotes}

\newcommand{\sqrtp}[2]{\sqrt[\leftroot{-2}\uproot{2}#2]{#1}}
\newcommand{\loss}{\mathcal{L}}
\newcommand{\encampsParant}[1]{\left(#1\right)}
\newcommand{\littleOh}[1]{o\encampsParant{#1}}
\newcommand{\bigOh}[1]{O\encampsParant{#1}}

\newtheorem{thm}{Theorem}
\newtheorem{lemma}{Lemma}

\begin{document}

\title{On Semi-Supervised Estimation of Distributions} 

\newcommand{\Dir}[1]{\text{Dir}_{#1}}
\newcommand{\maxmax}{\max} 
\newcommand{\minmin}{\min} 
\newcommand{\minmax}{min-max}
\newcommand{\data}{\{X_i\}^n_{i=1}}
\newcommand{\dataXY}{\{X_i, Y_i\}^n_{i=1}} 
\newcommand{\datap}{\{X'_i\}^m_{i=1}}
\newcommand{\datapXY}{\{X'_i, Y'_i\}^m_{i=1}}
\newcommand{\ind}[1]{ \bone\left\{#1\right\}}

\newcommand{\simiid}{\substack{ \text{} \\ \sim }}

\author{%
  \IEEEauthorblockN{Anonymous Authors}
}

%
 \author{%
   \IEEEauthorblockN{H.S.Melihcan Erol, Erixhen Sula and Lizhong Zheng}
   \IEEEauthorblockA{ Dept. EECS and RLE \\Massachusetts Institute of Technology\\
                     Cambridge, MA 02139\\
                     \{hsmerol, esula, lizhong\}@mit.edu}
                   }

\maketitle

\begin{abstract}
  We study the problem of estimating the joint probability mass function (pmf) over two random variables. In particular, the estimation is based on the observation of $m$ samples containing both variables and $n$ samples missing one fixed variable. We adopt the minimax framework with $l^p_p$ loss functions, and we show that the  composition of uni-variate minimax estimators achieves  minimax risk with the optimal first-order constant for $p \ge 2$, in the regime $m = o(n)$.
\end{abstract}
\section{Introduction}\label{sec:double-blind-policy}
Estimating the probability mass function (pmf) is a crucial statistical task. A commonly used formulation for estimation is  the minimax framework \cite{Wald49}. 
Early work resolved the minimax risk of pmf estimation under $l^2_2$ loss by identifying the minimax estimator \cite{trybula1958some,Olkin1979,wilczynski1985minimax}.
Later studies determined the constant of the first order for the minimax risk under KL-divergence, $l_1$ and $f$-divergences \cite{braess2004bernstein, HanJW14,pmlr-v40-Kamath15}. 

Meanwhile, the last decades have witnessed a rapid expansion in the sizes of the available datasets for which the labeling efforts lag behind. This results in heterogeneous datasets where a significant portion of the samples lack some of the variables. However,  the estimators considered in the works \cite{trybula1958some,Olkin1979,wilczynski1985minimax, braess2004bernstein, HanJW14,pmlr-v40-Kamath15}, can only operate in two modes in this setting: either the estimator should ignore the complete samples to estimate only the marginal pmf of the corresponding variable or the estimator should neglect the incomplete samples to form an estimate of the joint pmf. Adopting the naming conventions from machine learning literature \cite{pattern_recog}, these modes of operation can be categorized as unsupervised and supervised estimation, respectively. However, both modes suffer from inefficiencies, prompting the need for estimators that can leverage both labeled and unlabeled samples, commonly referred to as semi-supervised learning. This paper focuses on investigating the fundamental limits of semi-supervised pmf estimators.

In particular, we study the case where there are two random variables $X, Y$ jointly distributed with $p_{XY}$. We observe two datasets: $m$ i.i.d. samples of $(x_i, y_i)$ pairs drawn from $p_{XY}$ and $n$ samples of only $x_j$ drawn from the marginal distribution $p_X$.Our goal is to find the minimax estimator of $p_{XY}$ based on these observations. As a main result, we establish that the composition of minimax univariate estimators achieves the correct first-order term of the risk for the semi-supervised estimation problem in the regime $m= o(n)$.

The minimax pmf estimation problem with labeled and unlabeled samples remains unexplored in the existing literature. In the multivariate case, the analysis is complicated by nature's control over the number of samples with a fixed marginal. Previous works, such as \cite{He1990} and \cite{Kirichenko21}, have addressed related complications with slight variations. Unlike these works, where the number of samples is either generated from a fixed distribution or chosen adversarially, our study focuses on the case where the number of samples is generated from a distribution adversarially chosen by nature \cite{He1990, Kirichenko21}.
 
\section{  Preliminaries \& Notation }
 We employ Bachmann–Landau asymptotic notation and say $a_n = o(b_n)$ if $\lim\sup_n \frac{a_n}{b_n} = 0$, $a_n = O(b_n)$ if $\lim\sup_n \frac{a_n}{b_n} = K < \infty $,  $\Theta(a_n) = b_n$ if $a_n = O(b_n)$ and $b_n = O(a_n)$.   $\delta_x$ indicates a  unit point mass at $x \in \cX$.
We write for $x, y \in \mathbb{R}$, $x \wedge y \defeq \min(x,y)$, $x \vee y \defeq \max(x,y)$. We denote the space of probability distributions over the finite set $\cX$ by $\Delta_{\cX}$. We reserve the symbols $k_x = \abs{\cX}$ and $k_y = \abs{\cY}$. We use the upper case of a letter to denote a random variable and the lower case to indicate the realization of that random variable. 
 Let $\loss: \Delta_{\cX} \times \Delta_{\cX} \rightarrow  \mathbb{R}$ be a loss function. For a set of joint samples $s = \{(x_i, y_i)\}^m_{i=1}$, we define their marginal sets as $s_X \defeq \{x_i\}^m_{i=1}$ and $s_Y \defeq \{y_i\}^m_{i=1}$ and its conditional subsets as $s_{Y \mid X = x } = \{ (x_i, y_i) : 1 \le i \le m : x_i = x\}$. We denote by $T_x(s^m)$ the number of samples with value $x$ in the set $s^m$. We denote by $\hat{p}_X(s^m)$ the maximum likelihood estimator for the pmf estimation problem, which coincides with the empirical counts in $s^m$, i.e. $\hat{p}_X(x;s^m) = \frac{T_x}{m}$. 

In the minimax setting \cite{Wald49}, we assume that nature adversarially chooses a distribution $p_X$; $n$ samples $u^n = \{x_i\}^n_{i=1}$ are drawn i.i.d. from this distribution; and our goal is to design an estimator $\hat{q}: \cX^n \to \Delta_\cX$ based on the samples $u^n$ to minimize the expected loss. We denote the associated risk by $r^{\loss}_n$, and the problem is formulated as:  
\begin{equation}
    r^{\loss}_n \defeq \min_{\hat{q}_X} \max_{p_X \in \Delta_{\cX}} \Ed{U^n}{\loss(p_X, \hat{q}_X(U^n))} \label{eq:minimax_general}
\end{equation}
The choices of $l^2_2$, KL divergence, $l_1$ and general $f$-divergences for $\loss$ has been considered in the prior work \cite{trybula1958some,wilczynski1985minimax, braess2004bernstein, HanJW14,pmlr-v40-Kamath15}. In our case we chose $\loss$ to be the general $l^p_p$ norms for $p \ge 2$, namely for $p, q \in \Delta_{\cX}$:
\begin{equation}
\loss(p,q) = \|p - q\|^p_p \defeq \sum_{x \in \cX} \Big( \abs{ p(x) - q(x) } \Big)^p 
\end{equation}
and we denote the minimax estimation risk by $r^p_n$, then \eqref{eq:minimax_general} becomes:
\begin{equation}
r^p_n \defeq r^{l^p_p}_n = \min_{\hat{q}_X } \max_{p_X} \Ed{U^n}{ \|p_X - \hat{q}_X(U^n)\|^p_p} \label{eq:rmn2}
\end{equation}

The minimax risk scales as $r^p_n = \Theta(n^{-\frac{p}{2}})$. We defer the reader to the appendix of the full version of the paper for a proof 
This rate is achieved by the maximum likelihood estimator.  We will denote the constants of this rate by $\bar{C}_p \defeq \limsup_{n} n^{\frac{p}{2}} r^p_n $ and $\underline{C}_p = \liminf_n n^{\frac{p}{2}} r^p_n$. For the clarity of the presentation, we will group the parameters $\bar{C}_p$, $\underline{C}_p$ into $C_p$ and adopt the notation $\simeq$ to denote that $g_n \simeq C_p f_n$ if $ \limsup_n g_n / f_n \le \bar{C}_p$ and $\liminf_n g_n / f_n \ge \underline{C}_p$. Similarly, we say $g_n \lesssim C_p f_n$ if $\liminf_n g_n / f_n \le \underline{C}_p$ and $\limsup_n g_n / f_n \le \bar{C}_p$. 

This paper studies the extension of \eqref{eq:rmn2} in a multivariate semi-supervised setting. We assume that nature chooses a model $p_{XY}$ of a pair of jointly distributed random variables $X,Y$. The estimator has access to two datasets: a collection of complete $l^m \defeq \{(x'_i, y'_i)\}^{m}_{i=1}$ and a collection of incomplete samples  $u^n \defeq \{x_i\}^n_{i=1}$ generated i.i.d from $p_{XY}$ and $p_X$. Our goal is to design an estimator with minimal expected risk.  We formulate this new problem as: 
\begin{IEEEeqnarray}{rCl} \label{eq:Rmn}
R^p_{m,n} \defeq \min_{\hat{q}_{XY}}\max_{p_{XY} } \Ed{U^n, L^m}{ \|p_{XY} - \hat{q}_{XY}(U^n, L^m)  \|^p_p}  \IEEEeqnarraynumspace
\end{IEEEeqnarray}
Throughout our study of \eqref{eq:Rmn}, we will consider the following auxiliary problems:
\begin{gather*}
    R^p_m \defeq \min_{\hat{q}_{Y\mid X }}\max_{p_{XY}} \Ed{L^m}{\|p_{XY} - p_X \hat{q}_{Y \mid X }(L^m)\|^p_p}    \\
\Bar{R}^p_m \defeq \min_{\hat{q}_{Y \mid X}} \max_{p_{X}} \Ed{L^m_{X}}{ \max_{p_{Y\mid X}} \Ed{L^m_Y}{\|p_{XY} - p_X \hat{q}_{Y \mid X}(L^m)\|^p_p}} 
\end{gather*}
The problem $R^p_m$ corresponds to the limit of the problem $R^p_{m,n}$ as $n \rightarrow \infty$. Intuitively, in this case, there are sufficiently many incomplete samples to make the perfect estimation of $p_X$ possible. The difference between the two problems is that for $\bar{R}^p_m$, nature has an additional advantage in forming the dataset $l^m$: it can first observe the realization of the $X$ symbols, i.e. $l_X^m$, then choose $p_{Y|X}$ from which to generate the $Y$ values and finish the construction of the joint samples.

For a given univariate estimator $\hat{q}^*_n$, we define the conditional estimator $\hat{q}^{*, m}_{Y \mid X}( l^m )$ as the {\it conditional composition based on }$\hat{q}^*_n$ to be the concatenation of $\hat{q}^*_{n_x}(l_{Y|X=x})$, where $n_x$ is number of samples in $l_{Y\mid X =x}$, which is equal to $T_x(l^m_X)$. We further define the estimator $\hat{q}^{*, m, n}_{XY}$ as the {\it joint composition based on } $\hat{q}^*_n$ to be the estimator that estimates the conditional distribution $p_{Y\mid X}$ with the conditional composition estimator and the marginal $p_X$ with the ML estimator, i.e. $\hat{q}^{*,m,n}_{XY}(u^n, l^m) = \hat{p}_X(u^n \cup l^m_X)  \hat{q}^{m}_{Y \mid X}(l^m)$.

We say that an estimator  $\{\hat{q}_n\}$ is  \textit{first order minimax optimal} for the problem $r^p_n$ if $\max_{p_X}\Ed{U^n }{\|p_X - \hat{q}_X \|^p_p} = r^p_n + o(r^p_n)$.  The same definition carries over to the problems $R^p_{m,n}$,$R^p_m$, $\bar{R}^p_m$.

\section{Results}
\subsection{Main Theorems}
\begin{thm}  \label{thm:1} 
Let $\hat{q}^*_{n}$ be a minimax optimal estimator for $r^p_n$. Then the conditional composition $\hat{q}^{*,m}_{Y\mid X}$ based on $ \hat{q}^{*}_n$ is minimax optimal for $\bar{R}^p_m$:
\begin{IEEEeqnarray}{lCr}
\max_{p_{X}} \Ed{L^m}{\max_{p_{Y|X}}\|p_{XY} - p_X \hat{q}^{*,m}_{Y \mid X} \|^p_p} = \bar{R}^p_m 
\end{IEEEeqnarray}
\end{thm}
\begin{thm}\label{thm:2} Let $p\ge 2$ and  $\hat{q}^*_n$ be a first order minimax optimal estimator for $r^p_n$. Then the conditional composition $\hat{q}^{*,m}_{Y\mid X}$ based on $\hat{q}^{*}_n$ is first order minimax optimal for $R^p_m$:
\begin{IEEEeqnarray}{lCr}
 \max_{p_{XY}} \Ed{L^m}{ \|p_{XY} - p_X \hat{q}^{*,m}_{Y\mid X} \|^p_p }  = R^p_m + o\left( R^p_m \right)
\end{IEEEeqnarray}
\end{thm}
\begin{thm} \label{thm:3}  Let $m = o(n)$:
\begin{IEEEeqnarray}{rCl}
\abs{ R^p_{m,n} - R^p_m   } \le  O\left( {m^{-\frac{p -1 }{2}} (n)^{-1/2} } \right)
\end{IEEEeqnarray}
\end{thm}
\begin{thm}\label{thm:4}  Let $p \ge 2$ and $\hat{q}^*_n$ be a first-order optimal estimator for $r^p_n$. Then the joint composition $\hat{q}^{*,m,n}_{XY}$ based on $\hat{q}^{*}_n$ is first order minimax optimal for $R^p_{m,n}$ in the regime $m = o(n)$. 
\end{thm}
%
\subsection{Sketch of the proof}
The main result of this paper is given in \thmref{thm:4}. To establish \thmref{thm:4}, we first show in \thmref{thm:1} that the conditional composition of a minimax optimal estimator for $r^p_n$ is a minimax optimal estimator for $\bar{R}^p_m$. In \thmref{thm:2}, we connect the problems $\bar{R}^p_m$ and $R^p_m$. In particular, we show that when $p \ge 2$ the adversarial distribution of $\bar{R}^p_m$ is $\delta_x$, for which the problem $\bar{R}^p_m$ reduces to $R^p_m$. Finally in \thmref{thm:3} we show that $R^p_{m,n} = R^p_m + o(m^{-\frac{p}{2}})$ by studying the regime $m = o(n)$.
\section{Proofs for Theorems}
\subsection{Proof for \thmref{thm:1} }
We define:
\begin{IEEEeqnarray}{C}
    f(p_X, \hat{q}_{Y\mid X}) \defeq \Ed{L^m_X}{ \max_{p_{Y\mid X}} \Ed{L^m_Y}{\|p_{XY} - p_X \hat{q}_{Y\mid X}\|^p_p} }  \IEEEeqnarraynumspace \label{eq:defn_eq_fpx}
\end{IEEEeqnarray}
Let $\hat{q}^{**}_{Y\mid X} : (\cX \times \cY)^m \rightarrow (\Delta_{\cY})^{\abs{\cX}}$ be an  estimator for the conditional distribution $p_{Y\mid X}$. A sufficient condition for $\hat{q}^{**}_{Y\mid X}$ to achieve $\bar{R}^p_m$ is that for all $p_X$:
\begin{equation}
\min_{\hat{q}_{Y\mid X}} f(p_X, \hat{q}_{Y\mid X}) = f(p_X, \hat{q}^{**}_{Y \mid X}) \label{eq:suff}
\end{equation}
since
\begin{IEEEeqnarray*}{rCl}
\max_{p_X} f(p_X, \hat{q}^{**}_{Y \mid X})   &\ge& \bar{R}^p_m \defeq \min_{\hat{q}} \max_{p_X} \ f(p_X, \hat{q})  \\
\IEEEeqnarraymulticol{3}{r}{ \ge  \max_{p_X} \min_{\hat{q} } f(p_X, \hat{q}) =  \max_{p_X} f(p_X, \hat{q}^{**}_{Y\mid X}) }\IEEEeqnarraynumspace 
\end{IEEEeqnarray*}
where the first inequality is due to the substitution, the second inequality is the change of $\min \max$ with $\max \min $, and the equality follows from \eqref{eq:suff}. 
Now let us show \eqref{eq:suff} holds for the composition estimator $\hat{q}^{*,m}_{Y \mid X}$. Fix a $p_X$ and for the compactness of notation we define $p_{i,x} \defeq \prob{   T_x(L_X) = i}$. Then the left-hand side of \eqref{eq:suff} becomes:
\begin{IEEEeqnarray*}{rCl}
\IEEEeqnarraymulticol{3}{l}{= \min_{\hat{q}_{Y\mid X}}  \Ed{L^m_X}{\max_{p_{Y\mid X}} \Ed{L^m_Y}{\sum_{x \in \cX } \left( p_X(x) \right)^p \|p_{Y \mid X = x } -  \hat{q}_{Y\mid X=x }\|^p_p} }  }  \\
  \IEEEeqnarraymulticol{3}{l}{ =\sum_{x \in \cX } p^p_X(x)
 \min_{\hat{q}_{Y\mid X=x}}  \Ed{L^m_X}{\max_{p_{Y\mid X = x}} \Ed{L^m_Y}{\|p_{Y \mid X = x } -  \hat{q}_{Y\mid X=x }\|^p_p} }  } \IEEEeqnarraynumspace
 \IEEEyesnumber \label{eq:risk_substitute_pre} \\ 
&=& \sum_{x \in \cX }\left( p_X(x) \right)^p    \sum^m_{i=0}  p_{i,x}   r^p_i \IEEEyesnumber \label{eq:risk_substitute}
\end{IEEEeqnarray*}
In \eqref{eq:risk_substitute_pre} we note that  the optimization variables  are independent. Steps leading \eqref{eq:risk_substitute_pre} to \eqref{eq:risk_substitute} are given below and these steps demonstrate that \eqref{eq:suff} holds for the estimator $\hat{q}^{*,m}_{Y \mid X}$:
\begin{IEEEeqnarray*}{rCl}
    \IEEEeqnarraymulticol{3}{l}{  \min_{\hat{q}_{Y\mid X=x}}  \Ed{L^m_X}{\max_{p_{Y\mid X = x}} \Ed{L^m_Y}{\|p_{Y \mid X = x } -  \hat{q}_{Y\mid X=x }\|^p_p} } } \\
     &=&    \min_{\hat{q}_{Y\mid X=x}}  \sum^m_{i=0}   p_{i,x}  \max_{p_{Y\mid X = x}} \Ed{L^i_Y}{\|p_{Y\mid X =x } -\hat{q}^i_{Y\mid X =x }\|^p_p}   \\ 
     \IEEEeqnarraymulticol{3}{l}{
     = \sum^m_{i=0}  p_{i,x}  \min_{\hat{q}^i_{Y\mid X=x}} \max_{p_{Y\mid X = x}}  \Ed{L^i_Y}{\|p_{Y\mid X =x } -\hat{q}^i_{Y\mid X =x }\|^p_p}  } \IEEEeqnarraynumspace \IEEEyesnumber \label{eq:risk:another} \\
    \IEEEeqnarraymulticol{3}{l}{
     = \sum^m_{i=0}  p_{i,x}   \max_{p_{Y\mid X = x}}  \Ed{L^i_Y}{\|p_{Y\mid X =x } -\hat{q}^{*,i}_{Y\mid X =x }\|^p_p}  } \IEEEeqnarraynumspace \IEEEyesnumber \label{eq:risk:another_post}
\end{IEEEeqnarray*}
To establish \eqref{eq:risk:another_post}, we observe in \eqref{eq:risk:another} that the expression $\min_{\hat{q}^i_{Y\mid X=x}} \max_{p_{Y\mid X = x}}  \Ed{L^i_Y}{\|p_{Y\mid X =x } -\hat{q}_{Y\mid X =x }\|^p_p}$ is the problem $r^p_i$ and hence is achieved by $\hat{q}^{*,i}_{Y \mid X}$. 
\subsection{Proof for \thmref{thm:2} }
We start from \eqref{eq:risk_substitute} by recaliling the definition of $p_{i,x}$: 
\begin{IEEEeqnarray}{rCl}
\bar{R}^p_m &=& \max_{p_X} \sum_{x \in \cX} \sum^m_{i=0} {m \choose i} \left( p_X(x) \right)^{i+p} \left( 1 - p_X(x) \right)^{m-i} r^p_i \IEEEeqnarraynumspace \label{eq:use_prev} \\
&\simeq& \max_{p_X} \sum_{x \in \cX} C_p \left(   \frac{p_X(x)}{m} \right)^{\frac{p}{2}} + o(m^{-\frac{p}{2}})  \label{eq:use_conv} \\ 
&=& \frac{C_p}{m^{\frac{p}{2}}} + o(m^{-\frac{p}{2}}) \label{eq:subst_delta}
\end{IEEEeqnarray}
\eqref{eq:use_conv} follows from \lemref{lem:convergence_final}. In \eqref{eq:use_conv}, we note that for  $p \ge 2$ the problem is convex and symmetric in variables $\{ p_X(x) \}_{x \in \cX}$ in the first order. Therefore the optimizer is a vertex of the probability simplex, which leads to \eqref{eq:subst_delta}. We obtain the matching lower bound by substituting $p_X= \delta_{x}$ for some $x \in \cX$. We carry out the steps for completeness below:
\begin{IEEEeqnarray*}{rCl}
R^p_{m}
&\ge& \min_{\hat{q}_{Y \mid X }} \max_{p_{Y \mid X}} \Ed{L \sim \delta_x p_{Y \mid X}}{\|\delta_x p_{Y \mid X} - \delta_{x} \hat{q}_{Y \mid X}\|^p_p} \\ 
&=& \min_{\hat{q}_{Y \mid X }} \max_{p_{Y \mid X}} \Ed{L_Y \sim  p_{Y \mid X=x}}{\| p_{Y \mid X=x} -  \hat{q}_{Y \mid X= x}\|^p_p}  \\ 
&=& \min_{\hat{q}_{Y \mid X=x }} \max_{p_{Y \mid X=x}} \Ed{L_Y \sim  p_{Y \mid X=x}}{\| p_{Y \mid X=x} -  \hat{q}_{Y \mid X= x}\|^p_p}  = r^p_m
\end{IEEEeqnarray*}
Therefore by $R^p_m \le \bar{R}^p_m$ we obtain:
\begin{IEEEeqnarray*}{rCl}
 \frac{C_p}{m^{\frac{p}{2}}}  \lesssim   R^p_m \le \bar{R}^p_m \simeq \frac{C_p}{m^{\frac{p}{2}}} + o\left(\frac{1}{m^{\frac{p}{2}}}\right)
\end{IEEEeqnarray*}
finally we use \thmref{thm:1}.
\subsection{Proof for \thmref{thm:3} }
By \lemref{lem:thm3:upper}
$$
R^p_{m,n} - R^p_{m} \le  \gamma^p_{m,n}
$$
We note that  By \lemref{lem:Rmp_char} and \lemref{lem:rate_rnp} we have $r^p_{m+n} = \Theta{(n+ m)^{-\frac{p}{2}}}$ and $R^p_m = \Theta(m^{-\frac{p}{2}})$. Therefore In the regime $m=o(n)$, $\gamma^p_{m,n} = O(m^{-\frac{p-1}{2}} (n)^{-1/2})$. Finally, we observe that $R^p_{m,n}$ monotonically decreases $n$, and in the limit it reduces $R^p_m$.  An alternative proof for the lower bound $R^p_m \le R^p_{m,n}$ is given in \lemref{lem:thm3:lower}.
\subsection{Proof for \thmref{thm:4}}
By \thmref{thm:2}, the composition estimator $\hat{q}^{*,m}_{Y\mid X}$ is first order minimax optimal for $R^p_m$ when $p \ge 2$. Finally by \thmref{thm:3}, $R^p_{m,n}$ and $R^p_m$ has the same first order in the regime $m = o(n)$. 

\section{Supplementary Results}
For  \lemref{lem:convergence_final} and \lemref{lem:bernstein_trick} we introduce:
\begin{IEEEeqnarray}{rCl}
H^n_p(x) &\defeq& \sum^n_{i=0} {n \choose i } r^p_i x^{i+p} (1- x)^{n- i} \\
    G^n_p(x) &\defeq& \sum^n_{i=0} \binom{n}{i} r^p_i x^{i} \left(\frac{i}{n} \right)^p (1-x)^{n-i} 
\end{IEEEeqnarray}

\begin{lemma} \label{lem:convergence_final} 
$$ 
 H^n_{p}(x) = C_p \ \left( \frac{x}{n} \right)^{\frac{p}{2}} +  o(n^{-\frac{p}{2}}) 
$$
\end{lemma}
\begin{IEEEproof}
We fix the constant $c$ given in \lemref{lem:bernstein_trick}. 
There are two cases: 

In the first case   $x \ge c \frac{log^2(n)}{n}$:
\begin{IEEEeqnarray}{rCl}
H^n_{p}(x) &\simeq&  \sum^n_{i=0} \binom{n}{i}  \frac{C_p}{i^{\frac{p}{2}} + 1 }  x^{i+p} (1 -x )^{n -i}  \label{eq:lem1_1} \\
\IEEEeqnarraymulticol{3}{l}{= \frac{C_p}{n^{\frac{p}{2}}} \sum^n_{i=0} {n \choose i} \left(  \frac{i}{n} \right)^{\frac{p}{2}}  x^i (1- x)^{n -i}  + \bigOh{ \frac{ H^n_{p}(x)}{\sqrt{\log{n}}} } }  \IEEEeqnarraynumspace \label{eq:lem1_2} \\ 
\IEEEeqnarraymulticol{3}{l}{= \frac{C_p}{n^{\frac{p}{2}}} \left( x^{\frac{p}{2}} + O\left( n^{-1}  \right) \right) + \bigOh{\frac{ H^n_{p}(x)}{\sqrt{\log{n}}} }   } \label{eq:lem1_3} \\
\IEEEeqnarraymulticol{3}{l}{= C_p \left( \frac{x}{n}\right)^\frac{p}{2}  + \littleOh{n^{-\frac{p}{2}}}}\label{eq:lem1_4}
\end{IEEEeqnarray}
\eqref{eq:lem1_1} holds since $r^p_n \simeq {C_p} {n^{-\frac{p}{2}}}$ whereas in \eqref{eq:lem1_2} we use \lemref{lem:bernstein_trick}. To obtain \eqref{eq:lem1_3}, we utilize \label{thm:approximation_bernstein} as follows:  we set $f(x) = x^{\frac{p}{2}}$ and bound the error of $n$th order Bernstein polynomial approximation $B_n$ as: 
\begin{IEEEeqnarray*}{rCl}
    \abs{ B_n(x;f)- f(x)} &=& n^{-1} x (1-x) f''(x) /2  + o(n^{-1})  \\
    &=& n^{-1} p/4 (p/2 - 1) x^{\frac{p}{2} - 2 }    + o(n^{-1}) \\ 
    &\le&  n^{-1} p/4 (p/2-1) + o(n^{-1}) \IEEEyesnumber \label{eq:upper_bound} 
\end{IEEEeqnarray*}
where \eqref{eq:upper_bound} follows since $p \ge 2$. Therefore we conclude that convergence is uniform with error $O(n^{-1})$ for all $x \in (0,1)$. Finally, \eqref{eq:lem1_3} implies that $H^n_p(x) = O({n^{-\frac{p}{2}}})$ and in \eqref{eq:lem1_4}  we substitute this in the error term of \eqref{eq:lem1_3}.

For the second case we have  $x \le c  \frac{log^2(n) }{n}$: 
\begin{IEEEeqnarray}{rCl}
H^n_{p}(x) &=&   \sum^n_{i=0} {n \choose i} r^p_i x^{i + p} (1- x )^{n -i} \\
&\le&  C_p \ x^p  \sum^n_{i=0} {n \choose i} \frac{1}{i^{\frac{p}{2}} + 1} x^i (1 - x)^{n - i} \\ 
&\le& C_p\ x^p = \bigOh{ \frac{\log^{2p} (n) }{n^p} }  
\end{IEEEeqnarray}
Similarly $ C_p \left(\frac{x}{n}\right)^{\frac{p}{2}} = O\left( \frac{\log^{2p}(n)}{n^p}  \right) $ when $x \le c \frac{  \log^2{n} }{n}$, therefore $\abs{ C_p \left( \frac{x}{n} \right)^{\frac{p}{2}} - H^n_p(x) }  = \littleOh{n^{-\frac{p}{2}}} $. 
\end{IEEEproof}
\begin{lemma}\label{lem:bernstein_trick}
There exists a $c > 0$ such that for $x \ge  c  \frac{log^2 n }{n}$: 
\begin{equation}
\abs{H^n_p(x) - G^n_p(x)} =  \bigOh{ \frac{  H^n_p(x) }{ \sqrt{\log{n}} } }
\end{equation}
\begin{IEEEproof}
Fix $c > 0$, let $\delta_1,  \delta_2 > 0$, whose values will be determined later, we define $\Delta_p(x,i) \defeq \abs{x^p - \Big( \frac{i}{n} \Big)^{p}} $. As a result of triangular inequality:
\begin{IEEEeqnarray}{rCl}
  \abs{H^n_p(x) - G^n_p(x)} \le \sum^n_{i=0}  {n \choose i } x^i (1- x)^{n-i}  r^p_i\ \Delta_p(x,i) \label{eq:to_seperate} \IEEEeqnarraynumspace
\end{IEEEeqnarray}
Before analyzing this sum,  we note that by the mean value theorem, there exists $\xi \in (x \wedge \frac{i}{n}, x \vee \frac{i}{n} )$ hence $x^p - \left(\frac{i}{n} \right)^{p} = \left(x - \frac{i}{n}\right) \xi^{p-1}$, thus:
\begin{align} 
\Delta_p(x,i) = \Big|x^p - \Big(\frac{i}{n}\Big)^{p} \Big| \le p \Big|x -  \frac{i}{n} \Big|  \Big|x \vee \frac{i}{n} \Big|^{p-1} \label{eq:lem2:delta_upp}
\end{align}
Now we analyze the sum in \eqref{eq:to_seperate} over the ranges $ i \le nx$, $nx > i $ separately. 
 For the case $i \le nx$: 
    \begin{IEEEeqnarray}{rCl}
 &\ \ & \sum_{i\le nx }  {n \choose i } x^i (1- x)^{n-i}  r^p_i\ \Delta_p(x,i)    \\ 
    &=& \sum_{i \le n x - \delta_1} {n \choose i } x^i (1-x)^{n-i} r^p_i \Delta_p(x, i) \IEEEnonumber   \\ 
 \  &\ +&  \sum_{n x - \delta_1 < i < nx} {n \choose i } x^i (1-x)^{n-i} r^p_i \Delta_p(x, i) \label{eq:lem2:delta_subs} \\
 &\le& \sum_{i \le n x - \delta_1} {n \choose i } x^i (1-x)^{n-i} 2  \IEEEnonumber \\ 
 &\ +& \sum_{nx - \delta_1 \le i \le nx }  {n \choose i } x^i (1- x)^{n-i}  r^p_i \ p \abs{x - \frac{i}{n}} x^{p-1} \label{eq:lem2:delta_less_tail}
\IEEEeqnarraynumspace    \\ 
    &\le& 2 e^{-\frac{n}{x} \delta^2_1}  +
 \sum_{nx - \delta_1 \le i \le nx } p {n \choose i } x^i (1- x)^{n-i}  r^p_i \delta_1 x^{p-1}  \quad  \IEEEeqnarraynumspace \\ 
  &=& 2 e^{-\frac{n}{x} \delta^2_1}    + \frac{ p\ \delta_1}{x} H^n_{p}(x) \label{eq:bound_game}
   \end{IEEEeqnarray}
In \eqref{eq:lem2:delta_subs} we use \eqref{eq:lem2:delta_upp}.  In \eqref{eq:lem2:delta_less_tail}, we bound the lower tail of the binomial via \thmref{thm:binomial_tail} and observe that $\abs{x - \frac{i}{n}} \le \delta_1$ in the range $nx  -\delta_1 \le i \le nx$. We also note that in \eqref{eq:lem2:delta_less_tail}, \  $r^p_i \le 2$ and $\Delta_p(x,i) \le 1$.
Now we choose $\delta_1 = c_1 \sqrt{- \frac{x}{n} \log \left( \frac{1}{x} H^n_p(x)  \right) }$ and obtain:
\begin{IEEEeqnarray}{rCl} 
\eqref{eq:bound_game} &\le& 2 \ \left(  \frac{H^n_p(x)}{x} \right)^{c^2_1} + p\ c_1 \sqrt{- \frac{1}{n x} \log \left( \frac{1}{x} H^n_p(x)  \right) } H^n_p(x) \IEEEeqnarraynumspace    \IEEEnonumber
\end{IEEEeqnarray}
Hence to establish the lemma, we first show that $\sqrt{-\frac{1}{nx} \log{\frac{1}{x} H^n_p(x) }} = O( \frac{1}{  \sqrt{\log{n}} } )$. 
\begin{IEEEeqnarray*}{rCl}
\frac{1}{x} H^n_p(x) &\simeq& C_p \sum^n_{i=0} \frac{x^{p-1}}{i^{\frac{p}{2}} + 1} {n \choose i} x^i (1-x)^{n- i} \\
&\ge& C_p \frac{1}{n^{\frac{p}{2}} + 1 }  \left(\frac{ c \log^2{n}}{n} \right)^{^{p-1}} \sum^n_{i=0} {n \choose i} x^i (1-x)^{n-i}  \\ \IEEEyesnumber \label{eq:lower_bound}\\
&=& C_p \frac{1}{n^{\frac{p}{2}} + 1 }  \left(\frac{ c \log^2{n}}{n} \right)^{^{p-1}}\ge \frac{k'}{n^{\frac{3}{2}p}} \IEEEyesnumber \label{eq:Hpnx_lower_bound}
\end{IEEEeqnarray*}
In \eqref{eq:lower_bound}, we notice $x \ge \frac{c \log^2{n}}{n}$. On the right-hand side of \eqref{eq:Hpnx_lower_bound} we collect the constants in $k'$.
Therefore we establish that:
\begin{IEEEeqnarray}{rCl}
\sqrt{-\frac{1}{n x } \log{ \frac{1}{x} H^n_p(x) } } \le \sqrt{k'' \frac{1}{n x } \log{n} } \le  \bigOh{ \frac{1}{\log(n)} } \quad \IEEEeqnarraynumspace \label{eq:resulting_delta}
\end{IEEEeqnarray}
where in the first inequality we use that $x \ge c \frac{\log{n}}{n}$ and collect the constants in $k''$. Secondly, we need to show that $\frac{H^n_p(x)}{x}$ decays sufficiently fast. To this end, we have: 
\begin{IEEEeqnarray*}{rCl}
\IEEEeqnarraymulticol{3}{l}{ \sum_{i} \frac{C_p}{i^{\frac{p}{2}} + 1} {n \choose i} x^{i + p - 1} (1-x)^{n -i }  
\le \sum_{i} \frac{c_3}{i + 1} {n \choose i} x^i (1-x)^{n-i} } \\ 
&\le& c_3 \frac{ 1 - (1 - x)^{n+1}}{n+1} x^{p-2 } \le  c_3 \frac{x^{p-2}}{n+1}
\end{IEEEeqnarray*}
Therefore by choosing $c_1$ large enough we ensure that $B \left(\frac{H^n_p(x)}{x}  \right)^{c^2_1} = o(n^{-\frac{p}{2}})$.
Whereas in the  second case $i > n x$:
\begin{IEEEeqnarray}{rCl}
\IEEEeqnarraymulticol{3}{l}{  \sum_{nx < i \le nx +  \delta_2 } {n \choose i} x^i (1 -x )^{n-i} r^p_i \Delta_p (x, i)  } \IEEEnonumber  \\  
&\quad+& \sum_{  i >  nx + \delta_2  } {n \choose i} x^{i} (1-x)^{n-i} r^p_i \Delta_p(x,i)  \\ 
&\le& \sum_{nx < i \le nx +  \delta_2 } {n \choose i} x^i (1 -x )^{n-i} \abs{x - \frac{i}{n}} \left( \frac{i}{n} \right)^{p-1}   \IEEEnonumber \\  
&\quad+&2  \sum_{  i >  nx + \delta_2  } {n \choose i} x^{i} (1-x)^{n-i}  \\ 
&\le& \delta_2 \sum_{nx < i \le nx +  \delta_2 } {n \choose i} x^i (1 -x )^{n-i} \left( x + \delta_2 \right)^{p-1}    + 2 e^{-\frac{n \delta^2_2}{2 (x + \frac{\delta_2}{3} )}} \IEEEnonumber \\ 
\IEEEeqnarraymulticol{3}{l}{\le \delta_2 2^{p-1} \sum_{nx < i \le nx +  \delta_2 } {n \choose i} x^i (1 -x )^{n-i}   \left( x \vee \delta_2 \right)^{p-1} + e^{-\frac{n \delta^2_2}{\frac{2}{3} (x \vee \delta_2 )}}  }  \IEEEnonumber \\
\IEEEeqnarraymulticol{3}{l}{  \le  \delta_2\  2^{p-1} \sum_{nx < i \le nx +  \delta_2 } {n \choose i} x^i (1 -x )^{n-i}   \left( x \right)^{p-1}  + 2 e^{-\frac{n \delta^2_2}{\frac{2}{3} x }}  } \IEEEyesnumber
\label{eq:thm2:lower:last}
 \end{IEEEeqnarray}
Each step is justified in the corresponding step in the analysis for the range $i \le nx$, except now we are using the upper tail in \thmref{thm:binomial_tail}. In \eqref{eq:thm2:lower:last}, we see that the problem is identical to \eqref{eq:bound_game} except for the constants. Hence we choose $\delta_2 = c_2 \sqrt{-\frac{x}{n} \log \left( \frac{1}{x} H^n_p(x)  \right) }$ and  by \eqref{eq:Hpnx_lower_bound} we have $\delta_2 \le  c_3 \sqrt{ x \frac{\log{n} }{n} }$. This ensures that $\delta_2 \le x $ when $x \ge c \frac{\log^2{n}}{n}$ and the step \eqref{eq:thm2:lower:last} is valid.
\end{IEEEproof}
\end{lemma}

\begin{lemma} \label{lem:thm3:upper} 
For $p \ge 0$, there exists constants $\{c_i\}^p_{i=0}$ and $c', c''$ such that:
\begin{IEEEeqnarray*}{rCl}
R^p_{m,n} &\le& R^p_m + \gamma^p_{m,n} 
\end{IEEEeqnarray*}
with $\gamma^p_{m,n} = \sum^{\floor{p} - 1}_{i} c_i (R^p_m)^{\frac{p-i}{p}} \left( r^{p}_{m+n} \right)^{\frac{i}{p}}  
+ c' (R^p_m )^{\frac{p - \floor{p}}{p}} (r^p_{m+n})^{\frac{\floor{p}}{p}} 
+ c'' (R^p_m + r^p_{m+n})^{\frac{p - \floor{p}}{p}} (r^p_{n+m})^{\frac{\floor{p}}{p}} $
\\
\begin{IEEEproof}
First let us fix $\hat{q}_{XY}(U,L), U, L$ and let us define: 
\begin{equation}
\Gamma^{U, L}_{x,y}(u) \defeq p_{XY}(x,y) - u \hat{q}_{Y \mid X}(y \mid x)
\end{equation}
with derivatives:
\begin{IEEEeqnarray*}{rCl}
\abs{ \frac{d^i}{du^i} \abs{\Gamma^{U,L}_{x,y}(u)}^p} &=& p^{(i)} (\hat{q}_{Y\mid X}(y \mid x))^i \abs{\Gamma^{U,L}_{x,y}(u)}^{p-i} \\
&\le& p^{(i)} \abs{\Gamma^{U,L}_{x,y}(u)}^{p-i}  \label{eq:bound_derivatives}
\IEEEyesnumber
\end{IEEEeqnarray*}
\enlargethispage{-1cm}
Continuing, we have:
\begin{IEEEeqnarray}{rCl}
 &&\abs{\Gamma^{U,L}_{x,y}(\hat{q}_X(x))}^p \label{eq:above_lefthandside}\\
 &\le& \sum^{\floor{p} - 1 }_{i=0} \frac{p^{(i)}}{i!} \abs{\Gamma^{U,L}_{x,y}({p_X(x)})}^{p-i} \abs{\hat{q}_X(x) - p_X(x)}^i \IEEEeqnarraynumspace \IEEEnonumber  \\
&+& \Big( \abs{  \Gamma^{U,L}_{x,y}(p_X(x)) }^{p - \floor{p}} \vee \abs{ \Gamma^{U,L}_{x,y}(\hat{q}_X(x))}^{p - \floor{p}}  \Big)   \IEEEnonumber \\ && \abs{\hat{q}_X(x) - p_X(x) }^{ \floor{p} } \frac{p^{\floor{p}}}{\floor{p}!} \label{eq:bound_mvt} \\
&\le&\sum^{\floor{p} - 1 }_{i=0} \frac{p^{(i)}}{i!} \abs{\Gamma^{U,L}_{x,y}({p_X(x)})}^{p-i} \abs{\hat{q}_X(x) - p_X(x)}^i \IEEEeqnarraynumspace \IEEEnonumber \\
&+& \frac{p^{\floor{p}}}{\floor{p}!} \abs{ \Gamma^{U,L}_{x,y}(p_X(x))}^{p - \floor{p}} \abs{\hat{q}_X(x) - p_X(x) }^{\floor{p} }  \IEEEeqnarraynumspace \IEEEnonumber \\
&+& \frac{p^{\floor{p}}}{\floor{p}!} \abs{\Gamma^{U,L}_{x,y}(\hat{q}_X(x))}^{p - \floor{p}} \abs{\hat{q}_X(x) - p_X(x) }^{ \floor{p} }  \label{eq:above_righthandside}
\end{IEEEeqnarray}
\begin{align}
 &  \Ed{U,L}{\|p_{XY} - \hat{q}_{XY}\|^p_p}  \le \Ed{U,L}{\|p_{XY} - p_X \hat{q}_{Y\mid X}\|^p_p} \nonumber \\
&+\sum^{\floor{p} - 1 }_{i=1} \frac{p^{(i)}}{i!} \Ed{U,L}{ \sum_{x,y} \abs{\Gamma^{U,L}_{x,y}({p_X(x)})}^{p-i} \abs{\hat{q}_X(x) - p_X(x)}^i  }   \  \nonumber \\
&+ \frac{p^{\floor{p}}}{\floor{p}!} \Ed{U,L}{ \sum_{x,y} \abs{ \Gamma^{U,L}_{x,y}(p_X(x))}^{p - \floor{p}} \abs{\hat{q}_X(x) - p_X(x) }^{p - \floor{p} }   }  \   \nonumber \\
&+ \frac{p^{\floor{p}}}{\floor{p}} \Ed{U,L}{ \sum_{x,y} \abs{\Gamma^{U,L}_{x,y}(\hat{q}_X(x))}^{p - \floor{p}} \abs{\hat{q}_X(x) - p_X(x) }^{p - \floor{p} }   } \label{eq:sum_take_exp}  \\ 
&\Ed{U,L}{\|p_{XY} - \hat{q}_{XY}\|^p_p} \le \Ed{U,L}{\|p_{XY} - p_X \hat{q}_{Y\mid X}\|^p_p} \nonumber\\
&+\sum^{\floor{p} - 1 }_{i=1} \frac{p^{(i)} k^{\frac{i}{p}}_y}{i!} \E{ \| p_{XY} - p_X \hat{q}_{Y \mid X}\|^p_p  }^{\frac{p-i}{p}} \E{\|\hat{q}_X - p_X \|^p_p}^{\frac{i}{p}}   \ \nonumber  \\
&+ \frac{p^{\floor{p}} k^{\frac{\floor{p}}{p}}_y}{\floor{p}!}  \E{\|p_{XY} - p_X \hat{q}_{Y\mid X}\|^p_p}^{\frac{p- \floor{p}}{p}} \E{\|\hat{q}_X - p_X \|^p_p}^{\frac{\floor{p}}{p}}    \ \nonumber  \\
&+ \frac{p^{\floor{p}} k^{\frac{\floor{p}}{p}}_y}{\floor{p}!}  \E{\|p_{XY} - \hat{q}_X \hat{q}_{Y\mid X}\|^p_p}^{\frac{p - \floor{p}}{p}} \E{\|\hat{q}_X - p_X\|^p_p}^{\frac{\floor{p}}{p}} \label{eq:holder_step} \\\nonumber
&\defeq   \kappa_{m,n}(p_{XY}, \hat{q}_{XY}  )
\end{align}
In \eqref{eq:above_lefthandside} we Taylor expand the $h(u) \defeq \abs{\Gamma^{U,L}_{x,y}(u) (\hat{q}_{X}) }^p$ around the  $u = p_X(x)$  and use the upper bound \eqref{eq:bound_derivatives} for the derivatives. We bound the remainder of the Taylor expansion with the mean value theorem via the monotonicity of the derivatives as in \eqref{eq:lem2:delta_upp}, leading to \eqref{eq:bound_mvt}. In \eqref{eq:sum_take_exp} we sum over $x,y$ and take expectations with respect to $U,L \sim p_{XY}$ of both handsides \eqref{eq:above_lefthandside} ,\eqref{eq:above_righthandside} by noting that  $\Ed{U,L}{\sum_{x,y} \abs{\Gamma^{U,L}_{x,y} (\hat{q}(X) )}^p } = \Ed{U,L}{\|p_{XY} - \hat{q}_{XY} \|^p_p}$. To obtain \eqref{eq:holder_step}, we apply Hölder's inequality to each summation inside the expectations in \eqref{eq:sum_take_exp}. Taking maximum of both sides over $p_{XY}$ and taking the minimum of the left-hand side over $\hat{q}_{Y\mid X}$ we establish that for all $\hat{q}_{XY}$:
\begin{IEEEeqnarray*}{rCl}
\IEEEeqnarraymulticol{3}{l}{ R^p_{m,n} = \min_{\hat{q}_{XY}} \max_{p_{XY}} \Ed{U,L}{\|p_{XY}- \hat{q}_{XY}\|^p_p}  \IEEEeqnarraynumspace
} \\
&\le&  \max_{p_{XY}} \kappa_{m,n} (p_{XY}, \hat{q}_{XY})  \label{eq:pre_move_inside} \IEEEyesnumber
\end{IEEEeqnarray*}
\begin{IEEEeqnarray*}{rCl}
&\le& \max_{p_{XY}} \Ed{U,L}{\|p_{XY} - p_X \hat{q}_{Y\mid X}\|^p_p}\ +\\
\IEEEeqnarraymulticol{3}{l}{\sum^{\floor{p} - 1 }_{i=1} c_i \max_{p_{XY}} \E{ \| p_{XY} - p_X \hat{q}_{Y \mid X}\|^p_p  }^{\frac{p-i}{p}}
 \max_{p_{X}} \E{\|\hat{q}_X - p_X \|^p_p}^{\frac{i}{p}} } \\
&+& c'  \max_{p_{XY}}\E{\|p_{XY} - p_X \hat{q}_{Y\mid X}\|^p_p}^{\frac{p- \floor{p}}{p}} \max_{p_X} \E{\|\hat{q}_X - p_X \|^p_p}^{\frac{\floor{p}}{p}}  \\
&+& c''' \max_{p_{XY}} \E{\|p_{XY} - \hat{q}_X \hat{q}_{Y\mid X}\|^p_p}^{\frac{p - \floor{p}}{p}} \max_{p_X}\E{\|\hat{q}_X - p_X\|^p_p}^{\frac{\floor{p}}{p}} \\
\label{eq:move_max_inside} \IEEEyesnumber \\  
 R^p_{m,n}  &\le& \max_{p_{XY}} \Ed{U,L}{\|p_{XY} - p_X \hat{q}_{Y\mid X}\|^p_p}\  + \\
\IEEEeqnarraymulticol{3}{l}{ \sum^{\floor{p} - 1 }_{i=1} c_i \max_{p_{XY}} \E{ \| p_{XY} - p_X \hat{q}_{Y \mid X}\|^p_p  }^{\frac{p-i}{p}} 
 \max_{p_{X}} \E{\|\hat{q}_X - p_X \|^p_p}^{\frac{i}{p}} } \\
&+& c'  \max_{p_{XY}}\E{\|p_{XY} - p_X \hat{q}_{Y\mid X}\|^p_p}^{\frac{p- \floor{p}}{p}} \max_{p_X} \E{\|\hat{q}_X - p_X \|^p_p}^{\frac{\floor{p}}{p}}  \\
  &+&  c''  \max_{p_{XY}} \E{\|p_{XY} - p_X \hat{q}_{Y\mid X}\|^p_p}^{{\frac{p - \floor{p}}{p}}}  \max_{p_X}\E{\|\hat{q}_X - p_X\|^p_p}^{\frac{\floor{p}}{p}}   \\
  &+& c''  \max_{p_X} \E{\|p_X - \hat{q}_X\|^p_p}  \\
  \label{eq:the3_jensen_step} 
\IEEEyesnumber
\end{IEEEeqnarray*}
In \eqref{eq:move_max_inside}, we further upper bound \eqref{eq:pre_move_inside} by taking the maximum of each summation separately. We define the constants  are $c_i \defeq p^{(i)} k^{\frac{i}{p}}_y / i!, c'= c''' \defeq  p^{(\floor{p})} k^{\frac{\floor{p}}{p}}_y / \floor{p}!$ based on \eqref{eq:holder_step}. In \eqref{eq:the3_jensen_step} we note that:
\begin{gather*}
\max_{p_{XY}} \E{\|p_{XY} - \hat{q}_X \hat{q}_{Y\mid X}\|^p_p} \\  \le 2^{p-1} \max_{p_X} \E{ \| p_{XY} - p_X \hat{q}_{Y \mid X}\|^p_p  }  +  2^{p-1}  \max_{p_{XY}}\E{\|p_X - \hat{q}_X\|^p_p}
\end{gather*}
which is a consequence of convexity of  $\abs{x}^p $ for $p \ge 1$. For completeness, we include a proof for this in \lemref{lem:rate_rnp}. 
Finally, we choose $\hat{q}_X$, $\hat{q}_{Y\mid X}$ to be the minimax estimators of the $r^p_m$ and $R^p_m$ respectively to establish the \lemref{lem:thm3:upper}.
\end{IEEEproof}
\end{lemma}

\section{Conclusion}
In this work, we considered the problem of minimax pmf estimation problem under $l^p_p$ loss  when there are $m$ labeled and $n$ unlabeled samples. In particular, we showed that for $p\ge 2$,  the composition estimators of univariate minimax problems are optimal in the first order over the regime $m=o(n)$. Extending the results to $1 \le p \le 2$ and for $f$-divergences are possible future directions.  
\label{sec:conclusion}

\newpage
\bibliographystyle{IEEEtran}
\bibliography{references.bib}

\newpage \  \newpage
\appendix

\section{Supplementary Lemmas}
\begin{thm}[\cite{devore1993constructive},Theorem-3.1]  \label{thm:approximation_bernstein} If $f$ is bounded on A, differentiable in some neighborhood of $x$, and has second derivative $f^{\prime \prime}(x)$ for some $x \in A$, then
$$
\lim _{n \rightarrow \infty} n\left[B_n(f, x)-f(x)\right]=\frac{x(1-x)}{2} f^{\prime \prime}(x) .
$$
\end{thm}

\begin{thm}[\cite{Chung2002}, lemma-2] \label{thm:binomial_tail}
 Let $X_1, \ldots, X_n$ be independent random variables with
$$
\operatorname{Pr}\left(X_i=1\right)=p_i, \quad \operatorname{Pr}\left(X_i=0\right)=1-p_i .
$$
We consider the sum $X=\sum_{i=1}^n X_i$, with expectation $\mathrm{E}(X)=\sum_{i=1}^n p_i$. Then, we have:
\begin{align*}
& \operatorname{Pr}(X \leq \mathrm{E}(X)-\lambda) \leq e^{-\lambda^2 / 2 \mathrm{E}(X)} \\
& \operatorname{Pr}(X \geq \mathrm{E}(X)+\lambda) \leq e^{-\frac{\lambda^2}{2(\mathbb{E}(X)+\lambda / 3)}}
\end{align*}
\end{thm} 
\begin{lemma} \label{lem:Rmp_char}
We have
\begin{gather}
    \frac{C_p}{m^{\frac{p}{2}}} \le  R^p_m \le k_x\frac{C_p}{m^{\frac{p}{2}}}
\end{gather}
\begin{IEEEproof}
The lower bound follows from choosing $p_X = \delta_x$ and is considered in the proof for \thmref{thm:2}. The upper bound follows from: 
\begin{gather*}
\|p_X p_{Y\mid X=x} - p_X \hat{q}_{Y\mid X} \|^p_p =  \sum_{x \in \cX} ( p_X(x))^p \|p_{\mid X=x} - \hat{q}_{Y\mid X= x}\|^p_p 
\end{gather*}
yielding:
\begin{gather*}
\sup_{p_{XY}}\Ed{L^m}{\|p_X p_{Y\mid X=x} - p_X \hat{q}_{Y\mid X} \|^p_p} \\ \le  \sum_{x \in \cX} \sup_{p_X(x), p_{Y\mid X= x}} \left( p_X(x) \right)^p  \Ed{L^m}{\|p_{Y \mid X=x} - \hat{q}_{Y\mid X= x}\|^p_p } \\
\le  \sum_{x \in \cX} \sup_{ p_{Y\mid X= x}} \Ed{L^m}{ \|p_{Y \mid X=x} - \hat{q}_{Y\mid X= x}\|^p_p }
\end{gather*}
Finally taking $\min$ over $\hat{q}^{m}_{Y\mid X}$ and noting that variables $\{ \hat{q}^{*}_{Y\mid X=x} \}_{x \in \cX}$  are independent: 
\begin{gather*}
  \inf_{\hat{q}_{Y\mid X}}  \sup_{p_{XY}}\Ed{L^m}{\|p_X p_{Y\mid X=x} - p_X \hat{q}_{Y\mid X} \|^p_p} \\
  \le \sum_{x \in \cX} \inf_{\hat{q}_{Y\mid X=x}} \sup_{ p_{Y\mid X= x}} \Ed{L^m}{\|p_{Y \mid X=x} - \hat{q}_{Y\mid X= x}\|^p_p} 
  = k_x r^p_m
\end{gather*}
\end{IEEEproof}
\end{lemma}

\begin{lemma}\label{lem:rate_rnp}
There exists constants $c,C$ such that for all $n \ge 1$ and $p \ge 1$:
\begin{equation}
    \frac{c}{n^{\frac{p}{2}}} \le  r^p_n \le \frac{C }{n^{\frac{p}{2}}}
\end{equation}
\end{lemma}
\begin{IEEEproof}
From \lemref{lem:norm_equiv} we obtain:
\begin{gather}
\| p_X - \hat{q}_X(X^n_1) \|_p\ k_x^{ 1 - \frac{1}{p}} \ge   \| p_X - \hat{q}_X(X^n_1)\|_1  \\ 
  \| p_X - \hat{q}_X(X^n_1) \|_p \ge    \| p_X - \hat{q}_X(X^n_1)\|_1  \frac{1}{k_x^{1 - \frac{1}{p}}}  \\ 
  \| p_X - \hat{q}_X(X^n_1) \|^p_p \ge   \left(  \| p_X - \hat{q}_X(X^n_1)\|_1  \frac{1}{k_x^{1 - \frac{1}{p}}}  \right)^p  \\ 
  \| p_X - \hat{q}_X(X^n_1) \|^p_p \ge   \left(  \| p_X - \hat{q}_X(X^n_1)\|_1   \right)^p \frac{1}{k_x^{p-1}}   \\  
  \| p_X - \hat{q}_X(X^n_1) \|^p_p \ge   \left(  \| p_X - \hat{q}_X(X^n_1)\|_1   \right)^p k_x^{ 1- p }
\end{gather}
Now taking the expectation of both sides of the inequality over $X^n_1 \sim p_X$:
\begin{gather}
\E{\|p_X - \hat{q}_X(X^n_1)\|^p_p} \ge \E{ \left( \|p_X - \hat{q}_X(X^n_1)\| \right)^p } k_x^{ 1- p }  \nonumber \\
  \E{\|p_X - \hat{q}_X(X^n_1)\|^p_p} \ge \left( \E{\|p_X - \hat{q}_X(X^n_1)\|} \right)^p k_x^{ 1- p } \nonumber  \\ 
  \intertext{taking $\min\max$ of both sides: }
  \min_{ \hat{q}_X }\max_{p_X} \E{\|p_X - \hat{q}_X(X^n_1)\|^p_p}  \nonumber  \\  \ge  \min_{ \hat{q}_X }\max_{p_X} \left( \E{\|p_X - \hat{q}_X(X^n_1)\|} \right)^p k_x^{ 1- p } \label{eq:jensen_step}   \\ 
  \min_{ \hat{q}_X }\max_{p_X} \E{\|p_X - \hat{q}_X(X^n_1)\|^p_p} \nonumber \\ \ge  \left(  \min_{ \hat{q}_X }\max_{p_X}  \E{\|p_X - \hat{q}_X(X^n_1)\|} \right)^p k_x^{ 1- p }  \nonumber  \\ 
    \min_{ \hat{q}_X }\max_{p_X} \E{\|p_X - \hat{q}_X(X^n_1)\|^p_p}  \nonumber  \\ \ge  \left(  \sqrt{ \frac{2 (k-1)}{\pi n}} + O\left( \frac{1}{n^{3/4}} \right) \right)^p k_x^{ 1- p }  \label{eq:kamath_step} 
\end{gather} 
where in \eqref{eq:jensen_step}, we use Jensen's inequality, and we cite Corollary-9 from \cite{pmlr-v40-Kamath15} for the lower bound. Hence, 
\begin{align*}
   r^p_n \ge& \left(  \sqrt{ \frac{2 (k-1)}{\pi n}} + O\left( \frac{1}{n^{3/4}} \right) \right)^p k_x^{ 1- p } \\ 
   =&  \left( \left(  \sqrt{ \frac{2 (k-1)}{\pi n}}  \right)^p +  O\left( \frac{1}{n^{3/4}} \right)  \left( \sqrt{ \frac{2 (k-1)}{\pi n}} \right)^{p-1}  \right)    k_x^{ 1- p } \\ 
   =&  \left(  \sqrt{ \frac{2 (k-1)}{\pi n}}  \right)^p k_x^{ 1- p } +  O\left( \frac{1}{n^{ \frac{2p+ 1}{4}} } \right)  
\end{align*}
For the upper bound, we plug in the MLE estimator and use its moments:
\begin{align*}
    \prob{ \abs{ \hat{p}_X(x) - p_X(x) } \ge t  } \le 2e^{-2 n t^2}
\end{align*}
by Hoeffding's inequality, then for the MLE estimator $\hat{p}_X$: 
\begin{align*}
\E{ \|  \hat{p}_X - p_X \|^p_p } =& \sum_{x \in \cX } \E{   \abs{\hat{p}_X(x; X^n_1) - p_X(x)}^p } \\ 
=& \sum_{x \in \cX }  \E{  \abs{ \hat{p}_X(x; X^n_1) - p_X(x) }^p } \\ 
=& \sum_{x \in \cX }  \int^{\infty}_{t = 0} \prob{  \abs{ \hat{p}_X(x; X^n_1) - p_X(x) }^p  \ge t } dt  \\ 
=& \sum_{x \in \cX }  \int^{\infty}_{t = 0} \prob{  \abs{ \hat{p}_X(x; X^n_1) - p_X(x) }   \ge \sqrtp{t}{p} } dt \\
=& \sum_{x \in \cX }  \int^{\infty}_{u = 0} \prob{  \abs{ \hat{p}_X(x; X^n_1) - p_X(x) }   \ge u } p u^{p-1} du \\ 
=& \sum_{x \in \cX }  \int^{\infty}_{u = 0} 2 e^{-2 n u^2 } p u^{p-1} du \\ 
\intertext{We let $ v = 2 n u^2 $:  }  
=& \sum_{x \in \cX }  \int^{\infty}_{u = 0} 2 e^{-2 n u^2  } p \left(u \right)^{p-2} u du  \\ 
=& \sum_{x \in \cX }  \int^{\infty}_{u = 0} 2 e^{-v } p \left( \frac{v}{2n} \right)^{ \frac{p-2}{2}} \frac{1}{4n} dv \\ 
=&   \left( \frac{1}{2n} \right)^{\frac{p}{2}}  p \sum_{x \in \cX }  \int^{\infty}_{u = 0}  e^{-v }  \left( v\right)^{ \frac{p-2}{2}} dv \\ 
=&   \left( \frac{1}{2n} \right)^{\frac{p}{2}}  p \sum_{x \in \cX }  \Gamma\left( \frac{p}{2}\right) \\ 
=&  \left( \frac{1}{2n} \right)^{\frac{p}{2}}  p  k \Gamma\left(\frac{p}{2}\right)  \\ 
\le& \left( \frac{1}{2n} \right)^{\frac{p}{2}}  p  k \left( \frac{p}{2} \right)^{\frac{p}{2}}
\end{align*}
\end{IEEEproof}
\begin{lemma} \label{lem:norm_equiv} Let $\frac{p}{q} \ge 1$ and $\|\|_p, \|\|_q$ be the $l_p, l_q$ norms for $\reals^n$. Then $\forall x \in \reals^n$: 
\begin{align}
\|x \|_q \le n^{\frac{1}{q} - \frac{1}{p}} \| x \|_p 
\end{align}
\begin{IEEEproof}
This is a classical application of Hölder's inequality: 
\begin{align}
\|x\|^q_q =& \sum^n_{i=1} \abs{x_i}^q  1 \\
\le& \left(  \sum^n_{i=1}  \abs{x_i}^{p} \right)^{\frac{q}{p}} \left( \sum^{n}_{i=1} \abs{1}^{\frac{q}{p-1}} \right)^{1 - \frac{q}{p}}  \\ 
=& \left(  \sum^n_{i=1}  \abs{x_i}^{p} \right)^{\frac{q}{p}} n^{1 - \frac{q}{p}} 
\end{align}
Yielding: 
\begin{align}
   \|x\|_q  = \left( \sum^n_{i=1} \abs{x_i}^q \right)^{\frac {1}{q}} \le& \left( \left(  \sum^n_{i=1}  \abs{x_i}^{p} \right)^{\frac{q}{p}} n^{1 - \frac{q}{p}}  \right)^{\frac{1}{q}}  \\ 
   =&   \left(  \sum^n_{i=1}  \abs{x_i}^{p} \right)^{\frac{1}{p}} n^{\frac{1}{q} - \frac{1}{p}} \\
   =& \|x\|_p n^{\frac{1}{q} - \frac{1}{p}}
\end{align}
\end{IEEEproof}
\end{lemma}

\begin{lemma}\label{lem:thm3:lower} 
For all $m,n$ and $p \ge 1$:
$$
R^p_m \le R^p_{m,n}
$$
\begin{IEEEproof}
Let us denote the adversarial distribution of the $R^p_m$ by $p^*_{XY}$. We choose a prior $\Pi_{p_{XY}}$ over $\Delta_{\cX \times \cY}$ such that $\pi_{p_{XY}} = \delta_{p^*_X} \pi_{p_{Y \mid X}}$, we leave the choice of the $\pi_{p_{Y\mid X }}$ free as long as $S(\pi_{p_{Y\mid X = x }}) = \Delta_{\cY}$ for all $x \in \cX$. Then,
\begin{IEEEeqnarray}{rCl}
R^p_{m,n}&=&\min_{ \hat{q}_{XY}} \max_{p_{XY}} \E{\|p_{XY} - \hat{q}_{XY} \|^p_p } \\  
&\ge& \min_{\hat{q}_{XY}} \Ed{p_{XY} \sim \pi_{XY} }{ \Ed{U, L \sim p_{XY}}{\|\hat{q}_{XY} - p_{XY}\|^p_p}}  \label{eq:non_coop:bayes_step_1} \\ 
&=& \Ed{p_{XY} \sim \pi_{XY} }{ \Ed{U, L \sim p_{XY}}{\|\hat{q}^{\pi_{XY}}_{XY} - p_{XY}\|^p_p}} \label{eq:non_coop:bayes_step_2} \\
&=& \Ed{p_{XY} \sim \pi_{XY} }{ \Ed{U, L \sim p_{XY}}{\| p^*_X \hat{q}^{\pi_{XY}}_{Y \mid X} - p_{X} p_{Y\mid X}\|^p_p}} \label{eq:non_coop:bayes_step_3} \IEEEeqnarraynumspace \\ 
&=& \Ed{p_{XY} \sim \pi_{XY} }{ \Ed{U, L \sim p_{XY}}{\| p^*_X \hat{q}^{\pi_{XY}}_{Y \mid X} - p^*_{X} p_{Y\mid X}\|^p_p}} \label{eq:non_coop:bayes_step_4}  \\
&=& \min_{\hat{q}_{Y\mid X }} \max_{p_{Y\mid X}} \E{\|p^*_X \hat{q}_{Y \mid X } - p^*_X p_{Y\mid X } \|^p_p} \label{eq:non_coop:bayes_step_5} \\
&=& \min_{\hat{q}_{Y\mid X}} \max_{p_X p_{Y \mid X }} \E{\|p_X \hat{q}_{Y \mid X} - p_X p_{Y \mid X}\|^p_p} = R^p_m \label{eq:non_coop:bayes_step_6}
\end{IEEEeqnarray}
in \eqref{eq:non_coop:bayes_step_1} we lower bound the supremum with the average. In \eqref{eq:non_coop:bayes_step_2} we set
 $\hat{q}^{\pi_{XY}}_{XY}$ to be the Bayes estimator for the prior $\pi_{p_{XY}}$ which minimizes the posterior risk for an assignment $U, L$:
 \begin{IEEEeqnarray}{rCl}
\hat{q}^{\pi_{XY}}_{XY}(U,L) &\defeq& \argmin_{q_{XY} \in \Delta_{\cX \cY}} \Ed{p_{XY} \sim \pi_{ p_{XY} \mid U,L}}{\| p_{XY} - q_{XY} \|^p_p} \IEEEeqnarraynumspace  \\ 
&\defeq& \argmin_{q_{XY} \in \Delta_{\cX \cY} } F(q_{XY}, U, L)  
\end{IEEEeqnarray}
In  \eqref{eq:non_coop:bayes_step_3}, we note that  the functional $F(q_{XY}, U, L)$ is minimized by some element in the support of $\pi_{XY \mid U, L}$ by the convexity of $\|.\|^p_p$.  Since for all $U$ and $L$,  $S(  \pi_{XY \mid U, L} ) \subseteq \{p^*_X q_{Y \mid X} : q_{Y \mid X } \in \Delta_{\cY} \}$ we have $\hat{q}^{\pi_{XY}}_{X}(U, L) = p^*_X$. 
In \eqref{eq:non_coop:bayes_step_4} we again use that the $\pi_{X} = \delta_{p^*_X}$. In \eqref{eq:non_coop:bayes_step_5}, we use that the prior $\pi_{Y\mid X}$ is essentially free, and since any minimax risk can be approximated arbitrarily by a sequence of priors by the minimax theorem  \cite{Wald49}.  Finally \eqref{eq:non_coop:bayes_step_6} follows since $p^*$ is the adversarial distribution of $R^p_m$.  
\end{IEEEproof}
\end{lemma}

\begin{lemma}\label{lem:lp_triangle}
\begin{gather*}
\max_{p_{XY}} \E{\|p_{XY} - \hat{q}_X \hat{q}_{Y\mid X}\|^p_p} \\  \le 2^{p-1} \max_{p_X} \E{ \| p_{XY} - p_X \hat{q}_{Y \mid X}\|^p_p  }  +  2^{p-1}  \max_{p_{XY}}\E{\|p_X - \hat{q}_X\|^p_p}
\end{gather*}
\begin{IEEEproof}
We have,
\begin{gather}
\abs{p_{XY}(x,y) - \hat{q}_{XY}(x,y)}^p \label{eq:triangle_lp_LHS} 
\\ = \big| p_{XY}(x,y) - p_X(x) \hat{q}_{Y\mid X}(y \mid x) \nonumber \\
+ p_X(x) \hat{q}_{Y\mid X }(y \mid x ) - \hat{q}_{XY}(x,y) \big|^p  \nonumber \\
= 2^p \Big| \frac{1}{2} (  p_{XY}(x,y) - p_X(x) \hat{q}_{Y\mid X}(y \mid x) )  \nonumber  \\ + \frac{1}{2} ( p_X(x) \hat{q}_{Y\mid X }(y \mid x ) - \hat{q}_{XY}(x,y) ) \Big|^p  \nonumber \\
\le 2^{p-1}  \big|  p_{XY}(x,y) - p_X(x) \hat{q}_{Y\mid X}(y \mid x) \big|^p  \nonumber \\ \quad + 2^{p-1} \big| p_X(x) \hat{q}_{Y\mid X }(y \mid x ) - \hat{q}_{XY}(x,y) \big|^p \label{eq:jensen_step_triangle_lp}\\  
\le  2^{p-1}  \big|  p_{XY}(x,y) - p_X(x) \hat{q}_{Y\mid X}(y \mid x) \big|^p \nonumber  \\ + 2^{p-1} \big|p_X(x)  - \hat{q}_{X}(x) \big|^p  \label{eq:triangle_lp_RHS}
\end{gather}
where in \eqref{eq:jensen_step_triangle_lp} we use Jensen's inequality. Summing over $x,y \in \cX, \cY$ and taking expectation over $U,L \sim p_{XY}$ of both sides \eqref{eq:triangle_lp_LHS} and \eqref{eq:triangle_lp_RHS} we obtain: 
\begin{gather*}
\E{\|p_{XY} - \hat{q}_X \hat{q}_{Y\mid X}\|^p_p} \le \\
 2^{p-1}  \E{ \| p_{XY} - p_X \hat{q}_{Y \mid X}\|^p_p  }  +  2^{p-1} \E{\|p_X - \hat{q}_X\|^p_p}
\end{gather*}
taking maximum over both sides over $p_{XY}$ we obtain: 
\begin{gather*}
   \max_{p_{XY}} \E{\|p_{XY} - \hat{q}_X \hat{q}_{Y\mid X}\|^p_p} \le \\  \max_{p_{XY}} \left( 
 2^{p-1}  \E{ \| p_{XY} - p_X \hat{q}_{Y \mid X}\|^p_p  }  +  2^{p-1} \E{\|p_X - \hat{q}_X\|^p_p} \right)
\\
 \le 2^{p-1} \max_{p_{XY}} \E{ \| p_{XY} - p_X \hat{q}_{Y \mid X}\|^p_p  }  + 2^{p-1} \max_{p_X} \E{\|p_X - \hat{q}_X\|}
\end{gather*}     
\end{IEEEproof}
\end{lemma}

\end{document}